\documentclass[11pt]{article}
\usepackage[fleqn]{amsmath}
\usepackage[fleqn]{mathtools}
\usepackage[utf8]{inputenc}
\usepackage[english]{babel}
\usepackage{mathtools}
\usepackage{tabularx}
\usepackage{enumerate}
\usepackage{bbm}

\newcommand{\R}{\mathbb{R}}
\usepackage{amssymb,amsmath,amsfonts,,caption,setspace,hyperref}

\usepackage{float}
\usepackage{mathtools}
\usepackage{lipsum}
\usepackage{natbib}
\usepackage{breakcites}
\newtheorem{theorem}{Theorem}[section]
\newtheorem{lemma}{Lemma}[section]

\newtheorem{assumption}{Assumption}[section]
\newcommand{\proofend}{$\hfill\Box{~}$}
\newenvironment{Proof}{\noindent {\em{\bf Proof.}}}{\proofend\\}

\DeclarePairedDelimiter\norm{\lvert \lvert}{\rvert \rvert}
\newcommand{\argmin}[1]{\underset{#1}{\operatorname{arg}\,\operatorname{min}}\;}
\newcommand{\op}{\text{op}}

\usepackage{geometry}
\numberwithin{equation}{section}

\begin{document}
\begin{titlepage}
\title{\Large \bf High-dimensional inference robust to outliers with $\boldsymbol{\ell_1}$-norm penalization}
\author{
{\large Jad B\textsc{eyhum}}
\footnote{Corresponding author: jad.beyhum@gmail.com, Naamsestraat 69, 3000 Leuven, Belgium}\\\texttt{\small ORSTAT, KU Leuven}}
 
\date{}

\maketitle

\begin{abstract}
\noindent This paper studies inference in the high-dimensional linear regression model with outliers. Sparsity constraints are imposed on the vector of coefficients of the covariates. The number of outliers can grow with the sample size while their proportion goes to $0$. We propose a two-step procedure for inference on the coefficients of a fixed subset of regressors. The first step is a based on several square-root lasso $\ell_1$-norm penalized estimators, while the second step is the ordinary least squares estimator applied to a well-chosen regression. We establish asymptotic normality of the two-step estimator. The proposed procedure is efficient in the sense that it attains the semiparametric efficiency bound when applied to the model without outliers under homoscedasticity. This approach is also computationally advantageous, it amounts to solving a finite number of convex optimization programs.\\
\vspace{0in}\\
\noindent \textbf{KEYWORDS:}  robust regression, high-dimensional regression, $\ell1$-norm penalization, unknown variance.\\
\vspace{0in}\\
\noindent \textbf{MSC 2020 Subject Classification}: Primary 62J07 ; secondary 62J05, 62F35. 

\bigskip
\end{abstract}
\setcounter{page}{0}
\thispagestyle{empty}
\end{titlepage}
\pagebreak \newpage

\section{Introduction}
\label{sec:2}
The statistican observes a dataset of $n$  i.i.d. realizations of an outcome random variable $y_i$, a random vector of covariates $d_i$ with support in $\mathbb{R}^{K}$ and a random vector of controls $x_i$ with support in $\mathbb{R}^p$. $K$ is fixed while $p=p_n$ is allowed to go to infinity with the sample size. We assume that the following relationship holds:
 \begin{equation}\label{reg}y_i=d_i^{\top}\alpha+x_i^\top\beta+\gamma_i +u_i,\ \forall i=1,\dots,n,
\end{equation}
where $\alpha\in \mathbb{R}^K, \beta \in \mathbb{R}^{p}$, the error term $u_i$ is a real-valued random variable such that $\mathbb{E}[w_iu_i \vert \gamma_i=0]=0$ where $w_i=(d_i^\top, x_i^\top)^\top$ and $\gamma_i$ is a random variable. It also holds that $\{y_i,d_i,x_i,u_i,\gamma_i\}_i$ are i.i.d. and $\mathbb{E}[w_i w_i^\top\vert \gamma_i=0]$ exists and is positive definite. The observation $i$ is called an outlier if $\gamma_i\ne 0$. Let  $X=(x_1, \dots, x_n)^\top$, $\gamma = (\gamma_1,\dots,\gamma_n)^\top$ and $u=(u_1,\dots,u_n)$. The goal is to obtain inference results on the parameter $\alpha$ in the presence of outliers.

This model embodies many situations of practical interest. In Economics, $d_i$ could correspond to a binary treatment, while $x_i$ is a set of controls which need to be added to the regression to ensure exogeneity of $d_i$. The vector $\gamma_i$ may represent measurement errors, in this case, $\alpha$ is the average treatment effect of the whole population. Instead, if outliers arise because for some atypical individuals $d_i$ has a causal effect largely different from the one of the rest of the population, then $\alpha$ is the causal effect of the vast majority of the population. There may be many controls because the data is by nature high-dimensional or the statistician decides to include multiple transformations of a low number of variables to account for nonlinearities (see e.g. \cite{belloni2014inference}). In Biology, the researcher may guess that a few genes increase the risk of catching a given disease. To test this conjecture, a dataset of genomes of sick and healthy individuals is gathered. The variables in $d_i$ are the suspected genes while $x_i$ corresponds to the rest of the genome. Here, a likely cause for outliers is measurement errors.

This paper borrows from an approach to inference in the high-dimensional regression model which is rooted in the econometrics literature (see \citet{belloni2011?1, belloni2012sparse, belloni2014high, belloni2014inference, belloni2016inference, belloni2016post, belloni2017program, belloni2019valid}). We introduce the linear projections of each covariate in $d_i$ on the controls $x_i$. For $k=1,\dots, K$, let 
 \begin{equation}\label{lp}d_{ki}=x_i^\top\beta^k+\gamma_{i}^k +\xi_{i}^k,\ \forall i=1,\dots,n,
\end{equation}
where $ \beta^k \in \mathbb{R}^{p}$, the error term $\xi_{i}^k$ is a real-valued random variable such that $\mathbb{E}[x_i\xi_i^k \vert \gamma_i^k=0]=0$ and $\gamma_i^k$ is a random variable. All observations are i.i.d. and $\mathbb{E}[x_ix_i^\top\vert \gamma_i^k=0]$ exists and is positive definite. Again, the observation $i$ is called an outlier if $\gamma_i^k\ne 0$.  Let $\gamma^k = (\gamma_1^k,\dots,\gamma_n^k)^\top$ and $\xi^k=(\xi_1^k,\dots,\xi_n^k)^\top$. Similarly, the linear projection of $y_i$ on $x_i$ can be written  \begin{equation}\label{lpy}y_{i}=x_i^\top\beta^0+\gamma_{i}^0 +\xi_{i}^0,\ \forall i=1,\dots,n,
\end{equation}
where $\beta^0 = \sum_{k=1}^K\alpha_k\beta^k+\beta$, $\gamma^0 =\sum_{k=1}^K\alpha_k\gamma^k+\gamma$ and $\xi^0 =\sum_{k=1}^K\alpha_k\xi^k+u$.

We study a two-step estimation procedure. In the first step, we apply a variant of the square-root lasso estimator of \cite{belloni2011square} to the regressions in \eqref{lp} and \eqref{lpy}. The proposed variant has the advantage of being robust to outliers and allows to obtain estimates $\widehat{\xi}^0$ and $\widehat{\xi}^k$ of $\xi$ and $\xi^k$, respectively. We penalize the $\ell_1$-norm of both $\beta^k$ and $\gamma^k$. In the second step, we use the ordinary least squares estimator applied to the regression of $\widehat{\xi}^0$ on $\widehat{\xi}^1,\dots,\widehat{\xi}^K$. The rationale behind this second step lies in a moment condition satisfying the Neyman orthogonality condition (\cite{belloni2017program}). We show that, if the vectors $\beta^0,\dots, \beta^K$ are sparse enough and the proposition of outliers in all first-step regressions and in \eqref{reg} goes sufficiently quickly to $0$, the proposed two-step estimator of $\alpha$ is asymptotically normal which enables us to build tests and confidence intervals. Strikingly, the asymptotic variance of our estimator is the same as the one of the OLS estimator of $\alpha$ in the regression of $y_i$ on $d_i$ and $x_i$ (under the usual conditions of the linear regression model, in particular fixed $p$ and no outliers). In this sense the proposed estimator is efficient.

\textbf{Related literature.} This paper draws upon the literature in at least two different research field. The first is that of inference in the high-dimensional linear regression model, for which different approaches have been suggested. \citet{javanmard2014confidence, van2014asymptotically} and \cite{zhang2014confidence} propose to debias the LASSO estimator. \citet{ belloni2014inference} rely rather on a two-step approach similar to ours which attains the semiparametric efficiency bound. None of these approaches, however, has been shown to be robust to outliers.

The second related field is that of robust regression. Detailed accounts of this field can be found in \cite{rousseeuw2005robust,hampel2011robust} and \cite{maronna2018robust}. In the context of low-dimensional regression, the literature identifies a trade-off between efficiency and robustness, as explained below. $M$-estimators (such as the Ordinary Least-Squares (OLS) estimator) are often efficient when data are generated by the standard linear model with Gaussian errors and without outliers. However, this comes at the cost of robustness; $M$-estimators may be asymptotically biased in the presence of outliers. By contrast, $S$-estimators such as the Least Median of Squares (LMS) and the Least Trimmed Squares (LTS) are robust under several measures of robustness developed in the literature. They are also asymptotically normal in the model with Gaussian errors and without outliers but have a larger asymptotic variance than the OLS estimator in the standard linear model. In contrast, our estimation procedure yields an efficient estimator which can be asymptotically normal even in the presence of outliers.

Within the robust regression literature some authors have considered applying $\ell_1$-norm penalization to robust estimation. For low-dimensional linear regression (or nested special cases), see for instance \cite{gannaz2007robust, she2011outlier, lee2012regularization,lambert2011robust,dalalyan2012socp,li2012simultaneous,gao2016penalized, collier2017rate}. These works do not provide inference results. In this setup, \cite{beyhum2020inference} shows that a variant of the square-root lasso estimator is asymptotically normal and efficient. The present paper can be seen as an extension of this result in a high-dimensional context.

Recently, some authors have studied the problem of simultaneous estimation of the regression coefficients and the outliers when the number of variables can be larger than the sample size, see for instance \cite{alfons2013sparse, liu2017robust, virouleau2017high, yang2018general, liu2020high}. Closely related to our study is \cite{nguyen2012robust}, which proposed the extended lasso estimator which is a variant of the lasso estimator. \cite{dalalyan2019outlier} later refined their results. None of those works develop confidence intervals. Our main contribution is therefore to propose a $\ell_1$-norm penalized estimation procedure in the high-dimensional regression model which is robust to outliers, asymptotically normal and efficient. Moreover, in \cite{nguyen2012robust} the proposed theoretical choice of penalty level depends on the variance of the error term. Because our first step estimator is a variant of the square-root lasso estimator, we avoid this caveat and are able to propose a choice of penalty level which is not a function of the variance of the noise.

\textbf{Notation.} We use the following notations. For a matrix $M$, $M^{\top}$ is its transpose, $\norm{M}_2$,  $\norm{M}_1$ and $\norm{M}_{\infty}$ are the $\ell_2$-norm, $\ell_1$-norm and the sup-norm of the vectorization of $M$, respectively. $\norm{M}_{\text{op}}$ is the operator norm of $M$, $\norm{M}_0$ is the number of non-zero coefficients in $M$, that is its $\ell_0$-norm and $\norm{M}_{2,\infty}$ is the maximum of the $\ell_2$-norms of the columns of $M$. Moreover, $U(M)$ is the maximal diagonal element of $M$, $\lambda_{\min}(M)$ and $\lambda_{\max}(M)$ are its smallest and largest eigenvalue, respectively. When $M$ is a $n_1\times n_2$ matrix, $S\subset\{1,\dots,n_1\}$ and $T\subset\{1,\dots,n_2\}$, we denote by $M_{ST}$ the $|S|\times |T|$ submatrix of $M$ corresponding to the rows indexed by $S$ and the columns indexed by $T$. For a real number $x\in \mathbb{R}$, $\text{sign}(x)$ is equal to $1$ if $x\ge 0$ and $-1$ otherwise. For $a,b\in \R$, $a\vee b$ (resp. $a\wedge b$) denotes the maximum (resp. minimum) of $a$ and $b$. Next, for an integer $m$, $I_m$ is the identity matrix of size $m\times m$. Finally, $\mathbbm{1}_{\{\cdot\}}$ denotes the indicator function.
\section{Estimation}
 \subsection{Framework}
The probabilistic framework consists of a sequence of data generating processes (henceforth, DGPs) that depend on the sample size $n$. We consider an asymptotic setting where $n$ goes to $\infty$ while the number of $p=p_n$ of regressors is allowed (but does not have to) to go to infinity with $n$. The different regression coefficients $\beta,\beta^0,\dots,\beta^K$ and vectors $\gamma,\gamma^1,\dots,\gamma^K$ can vary with $n$, but $\alpha$ remains fixed.

The proposed estimation strategy is able to handle models where $\gamma^k$ is  sparse for all $k=1,\dots,K$, that is $\left|\left|\gamma ^k \right|\right|_0/n=o_P(1)$ or, in other words, $\epsilon^k\to 0$. Potentially, every outcome $y_i,d_i$ can be generated by a distribution that does not follow a linear model but the difference between the distribution of $y_i,d_i$ and the one yielded by a linear model can only be large for a negligible proportion of individuals. The subsequent theorems will help to quantify these statements.
 \subsection{Estimator}

We propose a two-step estimation strategy. In a first step, we estimate the regressions \ref{lp} and \ref{lpy}. For $k\in\{1,\dots,K\}$, $b\in\R^p$ and $c\in\R^n$, let 
$$Q^k(b,c)= \left\{\begin{array}{cc}  \frac{1}{n} \sum_{i=1}^n(y_i-x_i^\top b-c_i)^2& \text{if } k=0\\
\frac{1}{n} \sum_{i=1}^n(d_{ki}-x_i^\top b-c_i)^2& \text{otherwise}.
\end{array}\right.$$
We use the following estimators which have the advantages of being robust to outliers. 
\begin{equation}
\label{1step}
(\widehat{\beta}^k,\widehat{\gamma}^k)\in \argmin{b\in \mathbb{R}^{p},\ c\in \mathbb{R}^n} (Q^k(b,c))^{1/2} +\frac{\lambda_{\beta}^k}{n}\norm{\widehat{\Psi}b}_1 +\frac{\lambda_{\gamma}^k}{n} \norm{c}_1, \ \forall k=0,\dots,K,
\end{equation}
where $\{\lambda_\beta^k\}_{k=1}^K, \{\lambda_\gamma^k\}_{k=1}^K$ are sequences of positive penalty levels which properties will be specified below and $\widehat{\Psi}$ is the diagonal matrix with diagonal coefficients $\widehat{\Psi}_ {jj}=n^{-1/2}\sqrt{\sum_{i=1}^n x_{ji}^2}$.

The second step estimator is the ordinary least square estimator of the regression of $y_i - x_i^\top\widehat{\beta}^0-\widehat{\gamma}^0_i$ on $d_{1i}-x_i^\top\widehat{\beta}^1-\widehat{\gamma}^1_i, \dots,d_{Ki}-x_i^\top\widehat{\beta}^K-\widehat{\gamma}^K_i$, that is 
\begin{equation}\label{2step}
\widehat{\alpha}\in\argmin{a\in\R^K} \sum_{i=1}^n\left(\widehat{\xi}^0- \sum_{k=1}^Ka_k\widehat{\xi}^k_i\right)^2,
\end{equation}
where $\widehat{\xi}_i^0 = y_i - x_i^\top\widehat{\beta}^0-\widehat{\gamma}^0_i$ and for $k=1,\dots,K$, $\widehat{\xi}^k_i=d_{ki}-x_i^\top\widehat{\beta}^k-\widehat{\gamma}^k_i$. This estimator relies on the moment condition
\begin{equation}\label{moment} \mathbb{E}\left[(d_i-x_i^\top\beta^k-\gamma^k_i)\left((y_i-x_i^\top\beta^0-\gamma^0_i)-\sum_{k=1}^K\alpha_k(d_i-x_i^\top\beta^k-\gamma^k_i)\right)\right]=0,
\end{equation}
for all $k=1,\dots,K$. If $\mathbb{E}[x_iu_i]=0$ and $\mathbb{E}[x_i\xi^k_i]=0$ for $k=1,\dots,K$, the partial derivatives of the moment in \eqref{moment} with respect to $\beta^0,\dots,\beta^K$ and $\gamma^0,\dots,\gamma^K$ are $0$. As mentioned in the introduction, this moment therefore satisfies the Neyman orthogonality condition (see \cite{belloni2017program}). This property reduces the effect on the second step of mistakes made in the first step.

\subsection{Rate of convergence of the first step estimator} 
Let $k\in\{0,\dots,K\}$. In this subsection we derive the convergence rate of the estimator \eqref{1step} of the regression problem \eqref{lp} ($k\ne 0$) or \eqref{lpy} ($k=0$) . For a vector $v\in\R^p$ (resp. $v\in\R^n$), we denote by $v_T$ (resp. $v_S$), the vector such that $(v_T)_i= v_i$ (resp. $(v_S)_i=v_i$) if $\beta^k_i\ne 0$ (resp. $\gamma^k_i\ne0$) and $(v_T)_i= 0$ (resp. $(v_S)_i= 0$) otherwise. Moreover, we write $v_{T^c}= v-v_T$ (resp. $v_{S^c}= v-v_S$). Let $$\kappa^k = \min_{(h,f)\in\mathbb{C}^k} \frac{\frac{1}{\sqrt{n}}\norm{Xh+f}_2}{\norm{\widehat{\Psi}h}_2+\frac{1}{\sqrt{n}}\norm{f}_2}, $$ where
 $\mathbb{C}^k = \{(h,f)\in\mathbb{R}^p\times\R^n\ |\  \lambda_{\beta}^k\norm{\widehat{\Psi}h_{T^c}}_1+\lambda_{\gamma}^k\norm{f_{S^c}}_1 \le 3\lambda_{\beta}^k\norm{\widehat{\Psi}h_{T}}_1+3\lambda_{\gamma}^k\norm{f_{S}}_1\}.$  The constant $\kappa^k$ bears similarities with the extended restricted eigenvalue of \cite{nguyen2012robust} but has a different scaling. It is a generalization of the usual restricted eigenvalue (see \cite{bickel2009simultaneous}).

 Let $M^k=(\lambda_{\beta}^k\sqrt{\norm{\beta^k}_0})\vee( n\lambda_{\gamma}^k\sqrt{\epsilon^k})$. The following assumption is the key to derive rates of convergence of the estimator \eqref{1step}.
\begin{assumption}
\label{ascv}
The following holds:
\begin{enumerate}[\textup{(}i\textup{)}]  
 \item\label{cvii}we have $\lim\limits_{n\to\infty}\Pr\left(\lambda_{\beta}^k\ge 2\sqrt{n}\max\limits_{j=1,\dots,p}\frac{|(X^\top\xi^k)_j|}{\norm{\xi^k}_2\widehat{\Psi}_{jj}}\right)=1$ and $\lim\limits_{n\to\infty}\Pr\left(\lambda_{\gamma}^k\ge 2 \sqrt{n} \frac{\norm{\xi^0}_\infty}{\norm{\xi^0}_2}\right)=1$;
 \item\label{cvi} there exists $\kappa_*^k>0$ such that $\lim_{n\to \infty}\Pr(\kappa^k\ge\kappa_*^k)=1$;
 \item\label{cviii} $M^k/n=o(1)$;
\item\label{cviv} $\xi^k=O_{\textrm{P}}(\sqrt{n})$ .
 \end{enumerate}
\end{assumption}
Assumption \ref{ascv} \eqref{cvii} limits the choice of the penalty level. In practice, it is advisable to choose the lowest level of penalty satisfying this condition. \cite{belloni2011square} elaborates on how to choose $\lambda_{\beta}^k$ according to \eqref{cvii}. Lemma 2.3 and Corollary 2.4. in \cite{beyhum2020inference} provide guidance on how to pick  $\lambda_{\gamma^k}$ under this constraint. Condition \eqref{cvi} states that the extended restricted eigenvalue is bounded from below with probability approaching $1$. A sufficient condition when $x_i$ is a mean zero Gaussian vector is given in Lemma 1 in \cite{nguyen2012robust}, but this result does not allow the penalization of $\beta^k$ to depend on the regressors. Condition \eqref{cviii} is a joint sparsity constraint on $\beta^k$ and $\gamma^k$. Condition \eqref{cviv} is standard and holds if the entries of $\xi^k$ are i.i.d. with finite expectation because of the law of large numbers. Below, we give a sufficient condition for conditions \eqref{cvii} and \eqref{cvi} under a Gaussian design.

\begin{lemma} \label{sufficient}Assume that the following holds:
\begin{enumerate}[\textup{(}i\textup{)}]  
\item\label{si} $\{(x_i,\xi_i^k)\}_i$ are i.i.d and $x_i$ is independent of $\xi_i^k$;
 \item\label{sii} there exists a positive definite matrix $\Sigma$ such that $\lambda_{\max}(\Sigma)U(\Sigma)=O(1)$, $1/(\lambda_{\min}(\Sigma)\wedge \min_{k=1,\dots,p}\Sigma_{kk})=O(1)$ and $x_i$ are i.i.d. $\mathcal{N}(0,\Sigma)$;
 \item\label{siii} there exists $\sigma^k>0$ such that $\xi^k_i$ are i.i.d. $\mathcal{N}(0,(\sigma^k)^2)$;
 \item\label{siv} there exists $c>1$ such that $\lambda_{\beta}^k\ge2c\sqrt{n}\sqrt{\log(p)}$ and $\lambda_{\gamma}^k\ge 2c\sqrt{\log(n)}$;
\item\label{sv} $(\norm{\beta^k}_0\vee 1)\log(p)/n=o(1)$ and $\epsilon^k\log(n)=o(1)$.
 \end{enumerate}Then Assumption \ref{ascv} is satisfied.
\end{lemma}
In the particular setting of Lemma \ref{sufficient}, we see that condition \eqref{cviii} in Assumption \ref{ascv} becomes $\norm{\beta^k}_0=o(\sqrt{n}/\log(p))$ and $\epsilon^k=o(1/\log(n))$. As claimed in the introduction, the average proportion of outliers $\epsilon^k$ goes to $0$ while their number $n\epsilon^k$ can diverge. Let $\mu^k=(\norm{\beta^k}_0\vee \epsilon^k )$. The following theorem characterizes the rates of convergence of the first-step estimator.
\begin{theorem}\label{thcv}
Under Assumption \ref{ascv}, we have \begin{align*}
\norm{\widehat{\beta}^k-\beta^k}_1+\frac{1}{\sqrt{n}}\norm{\widehat{\gamma}^k-\gamma^k}_1&=O_{\textrm{P}}\left(\sqrt{\mu^k}\left(\frac{\norm{X}_{2,\infty}}{\sqrt{n}}\vee 1\right) \frac{M^k}{n}\right);\\
\norm{\widehat{\beta}^k-\beta^k}_2+\frac{1}{\sqrt{n}}\norm{\widehat{\gamma}^k-\gamma^k}_2&=O_{\textrm{P}}\left( \left(\frac{\norm{X}_{2,\infty}}{\sqrt{n}}\vee 1\right) \frac{M^k}{n}\right).
\end{align*}

\end{theorem}
The term $\norm{X}_{2,\infty}/\sqrt{n}$ is the maximum of the empirical standard deviations of the variables in $X$. It is equal to $\lambda_{\max}(\widehat{\Psi})$ and is an $O_{\textrm{P}}(1)$ when the variables in $X$ are uniformly bounded over $k$ and $n$ or when the assumptions of Lemma \ref{sufficient} hold (see Lemma \ref{chisquare} in the appendix).
Under the assumptions of Lemma \ref{sufficient}, we obtain the same rate of convergence of $\norm{\widehat{\beta}^k-\beta^k}_2+\frac{1}{\sqrt{n}}\norm{\widehat{\gamma}^k-\gamma^k}_2$ as in Corollary 1 of \cite{nguyen2012robust} (but without requiring knowledge of the variance of the error term), that is $$O_{\textrm{P}}\left(\sqrt{\frac{\norm{\beta^k}_0\log(p)}{n}}+\sqrt{\epsilon^k\log(n)}\right).$$

\subsection{second-step estimator}
In this section, we present sufficient assumptions for the asymptotic normality of $\widehat{\alpha}$ and consistent estimation of its asymptotic variance. Recall that $\widehat{\alpha}$ is the OLS estimator of the regression of $\widehat{\xi}^0$ on $\widehat{\xi}^1,\dots,\widehat{\xi}^K$. The first set of assumptions concerns the distribution of $(\xi^0,\dots,\xi^K)$ which we attempt to estimate in the first step. Let $\xi_i=(\xi_i^1,\dots,\xi_i^K)^\top$ and $\widehat{\xi}_i=(\widehat{\xi}_i^1,\dots,\widehat{\xi}_i^K)^\top$.
\begin{assumption}
\label{asan}
The following holds:
\begin{enumerate}[\textup{(}i\textup{)}]  
 \item\label{ani} $\{(d_i,x_i,u_i)\}_i$ are i.i.d. random variables;
 \item\label{anii} $\mathbb{E}[\xi_iu_i]=E[u_i]=0$;
 \item\label{aniii} $\Sigma_\xi = \mathbb{E}[\xi_i\xi_i^\top]$ exists and is positive definite;
\item\label{aniv} there exists $\sigma>0$ such that $\text{var}(u_i^2|\xi_i)=\sigma^2<\infty$. The conditional variance $\sigma^2$ does not scale with $n$.
 \end{enumerate}
\end{assumption}
These conditions are standard in the linear regression literature and guarantee that the OLS estimator of the regression of $\xi^0$ on $\xi^1,\dots,\xi^K$ is asymptotically normal.  Let us now introduce $\bar M = \max_{k=0}^K M^k$, $\bar\lambda_\beta =\max_{k=0}^K \lambda_\beta^k$, $\bar\lambda_\gamma =\max_{k=0}^K \lambda_\gamma^k$ and $\bar \mu=\max_{k=0}^K\mu^k$. The second set of assumptions ensures that $\xi^0,\dots,\xi^K$ are sufficiently well estimated in the first step.
\begin{assumption}\label{asrate}
The following holds:
\begin{enumerate}[\textup{(}i\textup{)}]  
 \item\label{asratei} $\bar M^2=o(n^{3/2})$;
 \item\label{asrateii} $\sqrt{\bar \mu} \left(\frac{\norm{X}_{2,\infty}}{\sqrt{n}}\vee 1\right)\bar M\left(\bar\lambda_\beta(\norm{X}_{2,\infty}/\sqrt{n}) + \bar\lambda_\gamma \sqrt{n}\right)=o_{\text{P}}(n^{3/2})$.
 \end{enumerate}
\end{assumption}
Under the assumptions of Lemma \ref{sufficient}, condition \eqref{asratei} in Assumption \ref{asrate} is satisfied if $\norm{\beta^k}_0=o(\sqrt{n}/\log(p))$ (the usual consistency condition of the lasso) and $\epsilon^k=o(1/(\sqrt{n}\log(n)))$. As already argued, we also have $(\norm{X}_{2,\infty}/\sqrt{n})=O_{\text{P}}(1)$ in this case. This implies that $$\sqrt{\bar \mu}\left(\frac{\norm{X}_{2,\infty}}{\sqrt{n}}\vee 1\right) \bar M\left(\bar\lambda_\beta(\norm{X}_{2,\infty}/\sqrt{n}) + \bar\lambda_\gamma\sqrt{n}\right)= O_{\text{P}}(\sqrt{\bar \mu} \bar M\left(\bar\lambda_\beta + \bar\lambda_\gamma \sqrt{n}\right))=O_{\text{P}}(\bar M^2)$$
and therefore Assumption \ref{asrate}\eqref{asrateii} is implied by Assumption \ref{asrate}\eqref{asratei}. These assumptions allow us to show the asymptotic normality of our two-step estimator. We have the following theorem.
\begin{theorem}\label{than}
Under assumptions \ref{ascv}, \ref{asan} and \ref{asrate}, we have 
$$\sqrt{n}(\widehat{\alpha}-\alpha)\xrightarrow{\Pr}\mathcal{N}(0,\sigma^2\Sigma_\xi^{-1}),\ \widehat{\sigma}\xrightarrow{\Pr}\sigma\  \text{and } \widehat{\Sigma}_\xi\xrightarrow{\Pr}\Sigma_\xi,$$
where $\widehat{\sigma}^2 = n^{-1}\sum_{i=1}^n(\widehat{\xi}^0_i-\sum_{k=1}^K\widehat{\alpha}_k\widehat{\xi}_i^k)^2$ and $\widehat{\Sigma}_\xi = n^{-1}\sum_{i=1}^n\widehat{\xi}_i\widehat{\xi}_i^\top$.
\end{theorem}
This result allows us to discuss efficiency. Let us consider the alternative problem of estimating $\alpha$, in the regression model \eqref{reg}, when there are no outliers ($\gamma,\gamma^1,\dots,\gamma^K=0$), $p$ is fixed and Assumption \ref{asan} holds. In this model, the OLS estimator of $(\alpha,\beta)^\top$ is asymptotically normal and the asymptotic variance of the OLS estimator of $\alpha$ (the projection on the first $K$ coordinates) is $\sigma^2\Sigma_\xi^{-1}$ by the Frisch-Waugh-Lovell theorem. Therefore, our estimator reaches the same asymptotic variance as the OLS estimator in this alternative model. In this sense, our estimator is efficient. This result is remarkable because it is obtained in a framework where there are outliers and $p$ can go to infinity.

An important remark concerns the meaning of confidence intervals developed using Theorem \ref{than}. They are obtained under an asymptotic setting with triangular array data in which the number of outliers is allowed to go to infinity while the proportion of outliers and nonzero coefficients of $\beta^0,\dots,\beta^K$ go to $0$. A 95\% confidence interval $I$ built with Theorem \ref{than} should be interpreted as follows: if the proportion of outliers in our data is low enough, the vectors $\beta^0,\dots,\beta^K$ are sparse and the sample size is large enough, then there is approximately a probability of $0.95$ that $\alpha$ belongs to $I$.

\section{Computation and simulations}
\subsection{Iterative algorithm}
\label{Iterative Algorithm}
We propose to use an algorithm similar to that of Section 5 in \cite{owen2007robust} to compute the first step estimators. Let $\widehat{\sigma}^k = Q^k(\widehat{\beta}^k,\widehat{\gamma}^k)$. Because $u=\min_{\sigma>0}\left\{\frac{\sigma}{2}+\frac{1}{2\sigma}u^2\right\}$, as long as $Q^k(\widehat{\beta}^k,\widehat{\gamma}^k)>0 $, we have
\begin{equation}\label{optiglob}
(\widehat{\beta}^k,\widehat{\gamma}^k, \widehat{\sigma}^k)\in \argmin{b\in \mathbb{R}^{p},c\in \mathbb{R}^n,s \in \mathbb{R}_+}\frac{s}{2}+ \frac{1}{2s}
Q^k(b,c)+\frac{\lambda^k_\beta}{n} \norm{\widehat{\Psi}b}_1+\frac{\lambda^k_\gamma}{n}\left|\left|c\right|\right|_1.
\end{equation}
This is a convex optimization program and the proposed approach is to iteratively minimize over $b$, $c$ and $s$. Let us start from $\left(b^{(0)},c^{(0)},s^{(0)}\right)$ and compute the following sequence for $t\in\mathbb{N}^*$ until convergence:
\begin{enumerate}
\item $b^{(t+1)}\in \argmin{b\in \mathbb{R}^{p}} Q^k(b,c^{(t)})+2\frac{\lambda^k_\beta s^{(t)}}{n} \norm{\widehat{\Psi}b}_1;$
\item $\gamma^{(t+1)}\in \argmin{c\in \mathbb{R}^n}  \left\vert\left| y-X\beta^{(t+1)}-c \right|\right\vert_2^2+\frac{2\lambda^k_{\gamma} s^{(t)}}{n}\left|\left|c\right|\right|_1;$
\item $s^{(t+1)}=\sqrt{Q^k(b,c)}$.
\end{enumerate}
Step 1 corresponds to the minimization of the usual lasso optimization program which is readily available in statistical softwares. The following lemma is a direct consequence of Section 4.2.2. in \cite{giraud2014introduction} and explains how to perform step 2.
\begin{lemma}\label{alphaiter}
For $i=1,\dots,n$, if $\left|y_i\mathbbm{1}_{k=0} + d_{ki}\mathbbm{1}_{k\ne 0} -x_i^\top b^{(t+1)} \right|\le\frac{\lambda^k_{\gamma} s^{(t)}}{n}$ then $c^{(t+1)}_i=0$, otherwise $c^{(t+1)}_i=y_i\mathbbm{1}_{k=0} + d_{ki}\mathbbm{1}_{k\ne 0} -x_i^\top b^{(t+1)}-\text{sign}\left(y_i\mathbbm{1}_{k=0} + d_{ki}\mathbbm{1}_{k\ne 0} -x_i^\top b^{(t+1)}\right)\frac{\lambda^k_{\gamma} s^{(t)}}{n}$.
\end{lemma}

\subsection{Simulations}
We apply our estimation procedure in a small simulation exercise. The $p$ regressors $x_i$ are i.i.d. $\mathcal{N}(0,I_p)$. The variable $d_i$ is unidimensional and generated as $d_i=x_i^\top\beta^1 +\gamma^1_i +\xi_i^1$, where $\beta^1\in\R^p$, $\beta^1_k=10$ for $6\le k\le 10$  and $\beta^1_k=0$ otherwise, $\xi^1_i$ are i.i.d. $\mathcal{N}(0,1)$ and $$\gamma^1_i=\left\{\begin{array}{cc} 0&\text{if $x_{11i}<\Phi^{-1}(1-\epsilon)$}\\
z&\text{if $x_{11i}\ge \Phi^{-1}(1-\epsilon)$},
\end{array}\right.$$
where $\Phi$ is the cumulative distribution function of the standard normal distribution and $z\in \R$. The outcome is given by $y_i=\alpha d_i + x_i^\top\beta +\gamma_i+\xi_i$, where $\beta\in\R^p$ is such that $\beta_k=10$ for $1\le k\le 5$ and $\beta_k=0$ otherwise, $\xi_i$ are i.i.d $\mathcal{N}(0,1)$ and $$\gamma_i=\left\{\begin{array}{cc} 0&\text{if $x_{6i}<\Phi^{-1}(1-\epsilon)$}\\
z&\text{if $x_{6i}\ge \Phi^{-1}(1-\epsilon)$}
\end{array}\right.$$
In tables ~\ref{fig:Cov} and ~\ref{fig:Cov1} we present the bias, the variance and the mean squared error (MSE) and the coverage of $95\%$ confidence intervals based on the asymptotic variance of Theorem \ref{than} for our estimator $\widehat{\alpha}$ for various values of $p,n,\epsilon$ and $z$. These quantities are computed as averages over 1,000 replications and in each replication the first-step estimators are computed using 10 iterations of the algorithm in Section \ref{Iterative Algorithm}.
The penalty levels are chosen according to Lemma \ref{sufficient}, that is $\lambda^k_\beta=2.02\sqrt{n}\sqrt{2\log(p)}$ and $\lambda^k_\gamma=2.02\sqrt{2\log(n)}$ for $k=0,1$. 
To compare our estimator to a procedure which is not robust to outliers, we report exactly the same information for the (biased) estimator $\widehat{\alpha}^{b}$ which is similar to $\widehat{\alpha}$ apart from the fact that it sets $\lambda^k_\gamma=0$ for $k=0,1$ ( $\lambda^k_\beta$ remains equal to $2.02\sqrt{n}\sqrt{2\log(p)}$). The confidence intervals for $\widehat{\alpha}^{b}$ are also computed using Theorem  \ref{than} but they are not asymptotically valid. We observe that our estimator has a very small prediction error and almost nominal coverage.
\begin{table}[!ht]
\begin{minipage}{0.48\textwidth}
  \centering
  \caption{$p=n=500$, $\epsilon=0.005$, $z=20$}
\label{fig:Cov} 
{\small
             \begin{tabular}{|c|c|c|}
  \hline
&$\widehat{\alpha}$ & $\widehat{\alpha}^{b}$ \\
  \hline
   bias  & $10^{-3}$&0.83\\
   var  & $8.10^{-3}$& 0.18\\
   MSE  &$8.10^{-3}$ & 0.86 \\
Coverage & 0.92 &  0.01\\
  \hline
      \end{tabular}  }
\end{minipage}
\begin{minipage}{0.48\textwidth}
  \centering
  \caption{$p=n=1,000$, $\epsilon=0.0025$, $z=40$}
\label{fig:Cov1} 
{\small
             \begin{tabular}{|c|c|c|}
  \hline
&$\widehat{\alpha}$ & $\widehat{\alpha}^{b}$ \\
  \hline
   bias  & $7.10^{-3}$&0.20\\
   var  & $3.10^{-3}$& 0.02\\
   MSE  &$3.10^{-3}$ & 0.06 \\
Coverage & 0.90 &  0.10\\
  \hline
      \end{tabular}  }
\end{minipage}
\end{table}


\section{Proofs}


\subsection{Proof of Lemma \ref{sufficient}}

\subsubsection{Proof that condition \eqref{cvii} in Assumption \ref{ascv} holds}
Conditional on $X$, $n^{-1/2}(X^\top\xi^k)_j=n^{-1/2}\sum_{i=1}^nx_{ij}\xi_i^k$ has a $\mathcal{N}(0, \widehat{\Psi}_{jj}^2(\sigma^k)^2)$ distribution. Therefore, by the Gaussian bound (see Lemma B.1 in \cite{giraud2014introduction}), we have 
\begin{equation}\label{gaussian_bound}\Pr\left(\frac{|(X^\top\xi^k)_j|}{\sqrt{n}\sigma^k\widehat{\Psi}_{jj}}\ge t\right)=E\left[\left. \Pr\left(\frac{|(X^\top\xi^k)_j|}{\sqrt{n}\sigma^k\widehat{\Psi}_{jj}}\ge t\right)\right| X\right]\le E[ e^{-\frac{t^2}{2}}|X]=e^{-\frac{t^2}{2}},\end{equation}
for $t\ge 0$.
Next, it holds that
\begin{align}
\notag &\Pr\left(\lambda_{\beta}^k< 2\sqrt{n}\max\limits_{j=1\dots,p}\frac{|(X^\top\xi^k)_j|}{\norm{\xi^k}_2\widehat{\Psi}_{jj}} \right)\\
\notag&\ge \Pr\left(c\sqrt{2\log(p)}<\max\limits_{j=1\dots,p}\frac{|(X^\top\xi^k)_j|}{\sqrt{n}\sigma^k\widehat{\Psi}_{jj}}+ \left|\frac{\sqrt{n}}{\norm{\xi^k}_2} - \frac{1}{\sigma^k}\right|\sup\limits_{j=1\dots,p}\frac{|(X^\top\xi^k)_j|}{\sqrt{n}\widehat{\Psi}_{jj}}\right)\\
\notag &\ge\Pr\left(\left(c-\frac{c-1}{2}\right)\sqrt{2\log(p)}< \max\limits_{j=1\dots,p}\frac{|(X^\top\xi^k)_j|}{\sqrt{n}\sigma^k\widehat{\Psi}_{jj}}\right)\\
 \label{pigeonhole} &\quad + \Pr\left(\frac{c-1}{2}\sqrt{2\log(p)}< \left|\frac{\sqrt{n}}{\norm{\xi^k}_2} - \frac{1}{\sigma^k} \right|\sup\limits_{j=1\dots,p}\frac{|(X^\top\xi^k)_j|}{\sqrt{n}\widehat{\Psi}_{jj}}\right),
\end{align}
by the pigeonhole principle. Now, because of \eqref{gaussian_bound}, we have 
\begin{align}\notag&\Pr\left(\left(c-\frac{c-1}{2}\right)\sqrt{2\log(p)} <\max\limits_{j=1\dots,p}\frac{|(X^\top\xi^k)_j|}{\sqrt{n}\sigma^k\widehat{\Psi}_{jj}}\right)\\
\notag &\quad=\mathbb{E}\left[\Pr\left(\left.\left(c-\frac{c-1}{2}\right)\sqrt{2\log(p)}< \max\limits_{j=1\dots,p}\frac{|(X^\top\xi^k)_j|}{\sqrt{n}\sigma^k\widehat{\Psi}_{jj}}\right| X\right)\right]\\
 \notag &\quad\le \mathbb{E}\left[\sum_{j=1}^{p}\Pr\left(\left.\left(c-\frac{c-1}{2}\right)\sqrt{2\log(p)}\le \frac{|(X^\top\xi^k)_j|}{\sqrt{n}\sigma^k\widehat{\Psi}_{jj}}\right| X\right)\right]\\
\label{limzero}&\le pe^{-\left(c-\frac{c-1}{2}\right)^2\log(p)}=e^{-\left(\left(c-\frac{c-1}{2}\right)^2-1\right)\log(p)}\to 0
\end{align}
because $c-\frac{c-1}{2}>1$.
Remark that by the law of large numbers and the continuous mapping theorem, $\sqrt{n}\norm{\xi^k}_2^{-1}\xrightarrow{\Pr} (\sigma^k)^{-1}$. Therefore, this and \eqref{limzero} yield
\begin{equation*}
\left|\frac{\sqrt{n}}{\norm{\xi^k}_2} - \frac{1}{\sigma^k} \right|\max\limits_{j=1\dots,p}\frac{|(X^\top\xi^k)_j|}{\sqrt{n}\widehat{\Psi}_{jj}}=o_{{\textrm{P}}}(1)O_{{\textrm{P}}}(\sqrt{\log(p)}),
\end{equation*}
which implies
\begin{equation}
\label{limzero2}\lim\limits_{n\to\infty}\Pr\left(\frac{c-1}{2}\sqrt{2\log(p)}< \left|\frac{\sqrt{n}}{\norm{\xi^k}_2} - \frac{1}{\sigma^k} \right|\max\limits_{j=1\dots,p}\frac{|(X^\top\xi^k)_j|}{\sqrt{n}\widehat{\Psi}_{jj}}\right)=0.
\end{equation}
By \eqref{pigeonhole}, \eqref{limzero} and \eqref{limzero2}, we obtain $\lim\limits_{n\to\infty}\Pr\left(\lambda_{\beta}^k\ge 2\max\limits_{j=1\dots,p}(|(X^\top\xi^k)_j|/(\norm{\xi^k}_2\widehat{\Psi}_{jj}))\right)=1$.

Next, we prove that $\lim\limits_{n\to\infty}\Pr\left(\lambda_{\gamma}^k\ge 2\sqrt{n}\norm{\xi^k}_\infty/\norm{\xi^k}_2\right)=1$. The Gaussian bound in Lemma B.1 of \cite{giraud2014introduction} yields 
\begin{equation}\label{gaussian_bound_2}\Pr\left(\frac{|\xi^k_j|}{\sigma^k}\ge t\right)\le e^{-\frac{t^2}{2}},\end{equation}
for $t\ge 0$.
Then, we have
\begin{align}
\notag &\Pr\left(\lambda_{\gamma}^k< 2\sqrt{n}\frac{\norm{\xi^k}_\infty}{\norm{\xi^k}_2}\right)\\
\notag&\ge \Pr\left(c\sqrt{2\log(n)}< \frac{\norm{\xi^k}_\infty}{\sigma^k}+ \left|\frac{\sqrt{n}}{\norm{\xi^k}_2} - \frac{1}{\sigma^k}\right|\norm{\xi^k}_\infty\right)\\
 \label{pigeonhole2}  &\ge\Pr\left(\left(c-\frac{c-1}{2}\right)\sqrt{2\log(n)}< \frac{\norm{\xi^k}_\infty}{\sigma^k}\right)
+ \Pr\left(\frac{c-1}{2}\sqrt{2\log(n)}<  \left|\frac{\sqrt{n}}{\norm{\xi^k}_2} - \frac{1}{\sigma^k}\right|\norm{\xi^k}_\infty\right),
\end{align}
by the pigeonhole principle. By \eqref{gaussian_bound_2}, we have $$\Pr\left(\left(c-\frac{c-1}{2}\right)\sqrt{2\log(n)}<\norm{\xi^k}_\infty/\sigma^k\right)\le e^{\left(\left(c-\frac{c-1}{2}\right)^2\right)\log(n)}\to 0.$$ This implies that $\norm{\xi^k}_\infty=O_{\textrm{P}}(\sqrt{\log(n)})$, which leads to $|\sqrt{n}\norm{\xi^k}_2^{-1}-(\sigma^k)^{-1}|\norm{\xi^k}_\infty=o_{\textrm{P}}(\sqrt{\log(n)})$ and shows that $\lim\limits_{n\to\infty}\Pr\left(\frac{c-1}{2}\sqrt{2\log(n)}<  \left|\sqrt{n}\norm{\xi^k}_2^{-1}-(\sigma^k)^{-1}\right|\norm{\xi^k}_\infty\right)=0$. We conclude using \eqref{pigeonhole2}.

\subsubsection{Proof that condition \eqref{cvi} in Assumption \ref{ascv} holds} 

In the rest of the proof, we denote $\norm{\beta^k}_0$  by $t$, $\norm{\gamma^k}_0$  by $s$ and $\lambda^k_\gamma/\lambda_\beta^k$ by $\lambda$. Let also $\Psi$ be the $p\times p$ diagonal matrix for which $\Psi_{kk}=\sqrt{\Sigma_{kk}}$ for $k=1,\dots,p$ and $C_X=\left(\sqrt{\lambda_{\min}(\Sigma)}/(8\norm{\Psi}_\infty)\right)\wedge 1/2$. The following proof bears similarities with that of Lemma 1 in \cite{nguyen2012robust} but is adapted to the specific penalty of our estimator. Take $(h,f)\in \mathbb{C}^k$. Let us first prove the following lemmas. 

\begin{lemma}\label{chisquare}Under the assumptions of Lemma \ref{sufficient}, it holds that $\lambda_{\max}(\widehat{\Psi})=O_{\textrm{P}}(1)$,  $1/\lambda_{\min}(\widehat{\Psi})=O_{\textrm{P}}(1)$ and \begin{equation} \label{boundmax}\lim\limits_{n\to\infty}\Pr\left(\frac{1}{\lambda_{\max}(\widehat{\Psi})}\ge \frac{3}{4\norm{\Psi}_\infty^2}\right).\end{equation}
\end{lemma}
\begin{Proof} Take $k\in\{1,\dots,p\}$. Because $X_i\sim\mathcal{N}(0,\Sigma)$, we have $\sum_{i=1}^nX_{ik}^2/\Sigma_{kk}\sim \chi^2_n$. Let us first prove $\lambda_{\max}(\widehat{\Psi})=O_{\textrm{P}}(1).$  Using Lemma 1 in \cite{laurent2000adaptive}, we obtain
$$
\Pr\left(\widehat{\Psi}_{kk}^2\ge \Sigma_{kk}+2\Sigma_{kk}\sqrt{\frac{\tau}{n}}+2\Sigma_{kk}\frac{\tau}{n}\right)\le \exp(-\tau).
$$
for $\tau>0$. Choose $\tau=2\log(p)$. Because $\lambda_{\max}(\widehat{\Psi})= \norm{\widehat{\Psi}}_\infty$, this implies 
\begin{align*} &\Pr\left(\lambda_{\max}(\widehat{\Psi})^2\ge \norm{\Psi}_\infty^2+2 \norm{\Psi}_\infty^2\sqrt{\frac{\tau}{n}}+2 \norm{\Psi}_\infty^2\frac{\tau}{n}\right)\\&\le \sum_{k=1}^p\Pr\left(\widehat{\Psi}_{kk}^2\ge \Sigma_{kk}+2\Sigma_{kk}\sqrt{\frac{\tau}{n}}+2\Sigma_{kk}\frac{\tau}{n}\right)\le p\exp(-\tau)\to 0.\end{align*}
Therefore, $\lambda_{\max}(\widehat{\Psi})=O_{\textrm{P}}(( \norm{\Psi}_\infty+2 \norm{\Psi}_\infty\sqrt{\tau/n}+2 \norm{\Psi}_\infty(\tau/n))^{1/2})$ which is an $O_{\textrm{P}}(1)$ by conditions \eqref{sii} and \eqref{sv} in Lemma \ref{sufficient}.

Next, we prove $1/\lambda_{\min}(\widehat{\Psi})=O_{\textrm{P}}(1)$. By Lemma 1 in \cite{laurent2000adaptive}, it holds that
\begin{equation} \label{chi2}
\Pr\left(\widehat{\Psi}_{kk}^2\le \Sigma_{kk}\left(1-2\sqrt{\frac{\tau}{n}}\right)\right)\le \exp(-\tau).
\end{equation}
We set $\tau =2\log(p)$. Because $1/\lambda_{\min}(\widehat{\Psi})=\max_{k=1,\dots,p}1/\widehat{\Psi}_{kk}$, we have 
$$\Pr\left(\frac{1}{\lambda_{\min}(\widehat{\Psi})^2}\ge \frac{1}{(1-2\sqrt{\frac{\tau}{n}})\min\limits_{k=1,\dots,p}\Sigma_{kk}}\right)\le \sum_{k=1}^p\Pr\left(\widehat{\Psi}_{kk}^2\le \Sigma_{kk}\left(1-2\sqrt{\frac{\tau}{n}}\right)\right)\le p\exp(-\tau)\to 0.$$
Remark that conditions \eqref{sii} and \eqref{sv} in Lemma \ref{sufficient} imply $1/((1-2\sqrt{\frac{\tau}{n}})\min\limits_{k=1,\dots,p}\Sigma_{kk}^2)=O_{\textrm{P}}(1)$, which shows that $1/\lambda_{\min}(\widehat{\Psi})=O_{\textrm{P}}(1)$.

Now, let us show \eqref{boundmax}. By \eqref{chi2}, we have $$\Pr\left(\frac{1}{\lambda_{\max}(\widehat{\Psi})}\ge  \frac{1}{(1-2\sqrt{\frac{\tau}{n}})\norm{\Psi}_\infty}\right)\le p\exp(-\tau).$$
Choosing $\tau=2\log(p)$, we obtain \eqref{boundmax} by condition \eqref{sv}.
\end{Proof}
\begin{lemma}\label{square}Under the assumptions of Lemma \ref{sufficient}, it holds that $$\lim\limits_{n\to\infty}\Pr\left(\norm{Xh}_2^2+\norm{f}_2^2 \ge n\frac{C_X^2}{2}\left(\norm{\widehat{\Psi}h}_2 + \frac{1}{\sqrt{n}}\norm{f}_2\right)^2\right)=1.$$ 
\end{lemma}
\begin{Proof} Take $v\in\R^p$, by Theorem 1 in \cite{raskutti2010restricted}, we have
\begin{align}\notag \frac{1}{\sqrt{n}}\norm{Xv}_2&\ge \frac{\sqrt{\lambda_{\min}(\Sigma)}}{4} \norm{v}_2-9\sqrt{U(\Sigma)}\sqrt{\frac{\log(p)}{n}}\norm{v}_1\\
\label{raskutti}&  \ge \frac{\sqrt{\lambda_{\min}(\Sigma)}}{4\lambda_{\max}(\widehat{\Psi})} \norm{\widehat{\Psi}v}_2-9\frac{\sqrt{U(\Sigma)}}{\lambda_{\max}(\widehat{\Psi})}\sqrt{\frac{\log(p)}{n}}\norm{\widehat{\Psi}v}_1.\end{align}
Next, because $(h,f)\in \mathbb{C}^k$, it holds that \begin{align*}\norm{\widehat{\Psi}h}_1&\le 4  \norm{\widehat{\Psi}h_T}_1+3\lambda  \norm{f_S}_1\\
&\le 4 \norm{h_T}_1+3\lambda  \norm{f_S}_1\\
&\le 4  \sqrt{t}\norm{\widehat{\Psi}h}_2+3\lambda  \sqrt{s}\norm{f}_2.
\end{align*} 
This and \eqref{raskutti} yield
\begin{align*}\frac{1}{\sqrt{n}}(\norm{Xh}_2+\norm{f}_2)&\ge \left(\frac{\sqrt{\lambda_{\max}(\Sigma)}}{4\lambda_{\max}(\widehat{\Psi})}-36\frac{\sqrt{U(\Sigma)}}{\lambda_{\max}(\widehat{\Psi})}\sqrt{\frac{t\log(p)}{n}} \right) \norm{\widehat{\Psi}h}_2\\ 
&\quad+\left(1-27\lambda\frac{\sqrt{U(\Sigma)}}{\lambda_{\max}(\widehat{\Psi})}\sqrt{s\log(p)}\right)\frac{1}{\sqrt{n}}\norm{f}_2.
\end{align*}
By conditions \eqref{sii}, \eqref{sv} and Lemma \ref{chisquare}, we have $\sqrt{U(\Sigma)t\log(p)}/(\sqrt{n}\lambda_{\max}(\widehat{\Psi}))=o_{\textrm{P}}(1)$ and $\lambda\sqrt{U(\Sigma)s\log(p)}/\lambda_{\max}(\widehat{\Psi})) =o_{\textrm{P}}(1)$, which implies
$$\lim\limits_{n\to\infty}\Pr\left(\norm{Xh}_2+\norm{f}_2 \ge \sqrt{n}C_X\left(\norm{\widehat{\Psi}h}_2 + \frac{1}{\sqrt{n}}\norm{f}_2\right)\right)=1.$$ 
We conclude the proof using the inequality $(a+b)^2\le 2(a^2+b^2)$.
\end{Proof}

We divide the set $\{1,\dots,p\}$ into subsets $T_1,\dots, T_q$ of size $t$, where $T_1=T$, $T_2$ contains the $t$ largest entries of $T^c$, $T_3$ contains the second $t$ largest entries of $T^c$ and so on. Let $s'\ge s$. In a similar manner, we split the set $\{1,\dots,n\}$ into $S_1,\dots,S_r$, where $S_1=S$, $S_2$ contains the $s'$ largest entries of $S^c$, $S_3$ contains the second $s'$ largest entries of $S^c$ and so on. We have 
\begin{align}
\notag\frac{1}{\sqrt{n}}|\left<Xh,f\right>|&\le \sum_{i,j}\frac{1}{\sqrt{n}}|\left<X_{S_iT_j}h_{T_j},f_{S_i}\right>|\\
\label{bounding}&\le \frac{1}{\sqrt{n}}\left(\max_{i,j}\norm{X_{S_iT_j}}_{\op}\right) \sum_{i=1}^q\norm{h_{T_i}}_2 \sum_{i=1}^r\norm{f_{S_i}}_2.
\end{align}
Let us show the following lemmas.
\begin{lemma}\label{wpa}Under the assumptions of Lemma \ref{sufficient}, it holds that
$$ \Pr\left(\max_{i,j} \frac{1}{\sqrt{n}}\norm{X_{S_iT_j}}_\op\le\sqrt{\lambda_{\max}(\Sigma)}\left(\sqrt{\frac{t}{n}}+ \sqrt{\frac{s'}{n}}+\tau\right)\right)\to 1$$
for any $\tau >0$.
\end{lemma}
\begin{Proof}
Remark that the rows of $X_{S_iT_j}$ have a $\mathcal{N}(0,\Sigma_{T_jT_j})$ distribution. Therefore, the entries of $X_{S_iT_j}\Sigma_{TjT_j}^{-1/2}$ are i.i.d. $\mathcal{N}(0,1)$. Let $\tau >0$. Applying Corollary 5.35 in \cite{vershynin2010introduction} to the matrix $X_{S_iT_j}\Sigma_{TjT_j}^{-1/2}$, we obtain that, with probability greater than $1-2\exp(-\tau^2n/2)$,
$\norm{X_{S_iT_j}\Sigma_{TjT_j}^{-\frac12}}_\op\le\sqrt{t}+ \sqrt{s'}+\tau\sqrt{n}$, which implies
\begin{equation}\label{opbound}
\Pr\left(\norm{X_{S_iT_j}}_\op\le\norm{\Sigma_{TjT_j}^{\frac12}}_{\op}\left(\sqrt{t}+ \sqrt{s'}+\tau\sqrt{n}\right)\right)\ge 1-2\exp(-\tau^2n/2)
\end{equation}
Taking the union bound over all $i$ and $j$, we have
$$ \Pr\left(\max_{i,j} \norm{X_{S_iT_j}}_\op\le\norm{\Sigma_{TjT_j}^{\frac12}}_{\op}\left(\sqrt{t}+ \sqrt{s'}+\tau\sqrt{n}\right)\right)\ge 1-2\binom{p}{t}\binom{n}{s'}\exp(-\tau^2n/2).$$
By condition \eqref{sv}, we have $\binom{p}{t}\le \left(\frac{ep}{t}\right)^t\le e^{t\left(\log\left(\frac{p}{t}\right)+1\right)}=o(e^n)$
and $\binom{n}{s'}\le\binom{n}{s} \le \left(\frac{en}{s}\right)^{s}\le e^{s\left(\log\left(\frac{n}{s}\right)+1\right)}=o(e^n).$
Therefore,  it holds that \begin{equation}\label{limitone}\Pr\left(\max_{i,j} \norm{X_{S_iT_j}}_\op\le\norm{\Sigma_{TjT_j}^{\frac12}}_{\op}\left(\sqrt{t}+ \sqrt{s'}+\tau\sqrt{n}\right)\right)\to 1.\end{equation}
Next, note that 
$$\norm{\Sigma_{TjT_j}^{\frac12}}_{\op}=\max\limits_{v\in\R^t,\norm{v}_2=1}\norm{\Sigma_{TjT_j}^{\frac12}v}_{2}\le \max\limits_{v\in\R^n,\norm{v_{T_j}}_2=1}\norm{\Sigma_{TjT_j}^{\frac12}v_{T_j}}_{2}
\le \sqrt{\lambda_{\max}(\Sigma)}. $$
This and \eqref{limitone} conclude the proof.
\end{Proof}
\begin{lemma}It holds that
$$ \sum_{i=1}^q \norm{h_{T_i}}_2\le  5\frac{\norm{\widehat{\Psi}h}_2 } {\lambda_{\min}(\widehat{\Psi})}+3\frac{\lambda}{\lambda_{\min}(\widehat{\Psi})}\sqrt{\frac{s'}{t}} \norm{f}_2\text{ and }\sum_{i=1}^r \norm{f_{S_i}}_2 \le 5\norm{f}_2 +\frac{3}{\lambda}\sqrt{\frac{t}{s'}}\norm{\widehat{\Psi}h}_2.$$
\end{lemma}
\begin{Proof}
We have
$\sum_{i=3}^q \norm{h_{T_i}}_2\le \sum_{i=3}^q \sqrt{t}\norm{h_{T_i}}_\infty 
\le\sum_{i=3}^q  \norm{h_{T_{i-1}}}_1/\sqrt{t}
\le \frac{1}{\sqrt{t}}\norm{h_{T^c}}_1/\sqrt{t}
$.
Because $(h,f)\in \mathbb{C}^k$, it holds that $\norm{\widehat{\Psi}h_{T^c}}_1\le 3\sqrt{t}\norm{\widehat{\Psi}h}_2+ 3\lambda\sqrt{s'}\norm{f}_2$.
This yields
\begin{align*}
\sum_{i=1}^q \norm{h_{T_i}}_2&\le 2\norm{h}_2 +\sum_{i=3}^q \norm{h_{T_i}}_2\\
&\le 2\norm{h}_2 +\frac{\norm{h_{T^c}}_1}{\sqrt{t}}\\
&\le  2\norm{h}_2 +\frac{\norm{\widehat{\Psi}h_{T^c}}_1}{\lambda_{\min}(\widehat{\Psi})\sqrt{t}}\\
&\le 2\norm{h}_2 +\frac{3}{\lambda_{\min}(\widehat{\Psi})\sqrt{t}}\left(\sqrt{t}\norm{\widehat{\Psi}h}_2 + \lambda \sqrt{s}\norm{f}_2\right)\\
&\le 5\frac{\norm{\widehat{\Psi}h}_2 } {\lambda_{\min}(\widehat{\Psi})}+3\frac{\lambda}{\lambda_{\min}(\widehat{\Psi})}\sqrt{\frac{s'}{t}} \norm{f}_2.
\end{align*}
Similarly, we have $\sum_{i=3}^r \norm{f_{T_i}}_2\le \sqrt{s'} \norm{f_{T^c}}_1$ and $\norm{f_{T^c}}_1\le (3/\lambda)\sqrt{t}\norm{\widehat{\Psi}h}_2+ 3\sqrt{s'}\norm{f}_2$, which implies
\begin{align*}
\sum_{i=1}^r \norm{f_{S_i}}_2&\le 2\norm{f}_2 +\sum_{i=3}^r \norm{f_{S_i}}_2\\
&\le 2\norm{f}_2 +\frac{\norm{f_{T^c}}_1}{\sqrt{s'}}\\
&\le 2\norm{f}_2 +\frac{3}{\sqrt{s'}}\left(\frac{\sqrt{t}}{\lambda}\norm{\widehat{\Psi}h}_2+ \sqrt{s}\norm{f}_2\right)\\
&\le 5\norm{f}_2 +\frac{3}{\lambda}\sqrt{\frac{t}{s'}}\norm{\widehat{\Psi}h}_2.
\end{align*}
\end{Proof}
\begin{lemma}\label{scalar}Under the assumptions of Lemma \ref{sufficient}, it holds that
$$ \Pr\left(|\left<Xh,f\right>|\le n\frac{C_X^2}{8}\left(\norm{\widehat{\Psi}h}_2+\frac{1}{\sqrt{n}}\norm{f}_2\right)^2\right)\to 1.$$
\end{lemma}
\begin{Proof}
For $\tau >0$, let us work of the event 
$$ \max_{i,j} \frac{1}{\sqrt{n}}\norm{X_{S_iT_j}}_\op\le\sqrt{\lambda_{\max}(\Sigma)}\left(\sqrt{\frac{t}{n}}+ \sqrt{\frac{s'}{n}}+\tau\right)$$
which proability goes to $1$ with $n$ according to Lemma \ref{wpa}. By \eqref{bounding}, $n^{-1/2}|\left<Xh,f\right>|$ is upper bounded by 
\begin{align*}
&\sqrt{\lambda_{\max}(\Sigma)}\left(\sqrt{\frac{t}{n}}+ \sqrt{\frac{s'}{n}}+\tau\right)\left(5\frac{\norm{\widehat{\Psi}h}_2 } {\lambda_{\min}(\widehat{\Psi})}+3\frac{\lambda}{\lambda_{\min}(\widehat{\Psi})}\sqrt{\frac{s'}{t}} \norm{f}_2 \right)\left(5\norm{f}_2 +\frac{3}{\lambda}\sqrt{\frac{t}{s'}}\norm{\widehat{\Psi}h}_2\right)\\
&\le 25 \sqrt{n}\frac{\sqrt{\lambda_{\max}(\Sigma)}}{\lambda_{\min}(\widehat{\Psi})}\left(\sqrt{\frac{t}{n}}+ \sqrt{\frac{s'}{n}}+\tau\right)  \max\left(1, \lambda\sqrt{\frac{s'}{t}},\frac{1}{\sqrt{n}\lambda}\sqrt{\frac{t}{s'}}\right)\left(\norm{\widehat{\Psi}h}_2+\frac{1}{\sqrt{n}}\norm{f}_2\right)^2.
\end{align*}
Take $s'=\max(s,\log(p) t)$. By condition \eqref{sv} and Lemma \ref{chisquare}, we have 
$$ \frac{1}{\lambda_{\min}(\widehat{\Psi})}\max\left(\sqrt{\frac{t}{n}},\sqrt{\frac{s'}{n}}, \lambda\sqrt{\frac{s'}{t}},\frac{1}{\sqrt{n}\lambda}\sqrt{\frac{t}{s'}}\right) =o_{\textrm{P}}(1).$$
Therefore, choosing $\tau$ sufficiently small, we obtain the result.
\end{Proof}

\noindent Let us now conclude the proof of Lemma \ref{sufficient}. It holds that 
$$\norm{Xh+f}_2^2= \norm{Xh}^2+\norm{f}_2^2+ 2\left<Xh,f\right>.$$
Therefore, by lemmas \ref{square} and \ref{scalar}, we have
$$\lim\limits_{n\to 1\infty}\Pr\left( \frac{1}{n}\norm{Xh+f}_2^2\ge \frac{C_X^2}{4}\left(\norm{\widehat{\Psi}h}_2+\frac{1}{\sqrt{n}}\norm{f}_2\right)^2\right)=1,$$
which yields the result with $\kappa_*^k=C_X/2$.

\subsection{Proof of Theorem \ref{thcv}}

We only prove the result for $k=0$. The proof for $k\in\{1,\dots,K\}$ being similar. Throughout this proof, we work on the event 
$$\left\{\kappa^0>\kappa_*^0\right\}\cap \left\{\lambda_{\beta}^0\ge 2\sqrt{n}\sup\limits_{j=1\dots,p}\frac{|(X^\top\xi^0)_j|}{\norm{\xi^0}_2\widehat{\Psi}_{jj}}\right\}\cap\left\{\lambda_{\gamma}^0\ge 2\sqrt{n}\frac{\norm{\xi^0}_\infty}{\norm{\xi^0}_2}\right\}\cap \left\{\kappa_*^0>2\frac{M^0}{n}\right\} ,$$ 
which probability goes to $1$ because of Assumption \ref{ascv}. Let us define $h=\widehat{\beta}^0 -\beta^0$ and $f=\widehat{\gamma}^0-\gamma^0$. Now, remark that
\begin{align}
 \notag \norm{\widehat{\Psi}\beta^0}_1- \norm{\widehat{\Psi}\widehat{\beta}^0}_1&= \norm{\widehat{\Psi}\beta^0}_1-\norm{\widehat{\Psi}\beta^0+\widehat{\Psi}h}_1\\
\notag&= \norm{\widehat{\Psi}\beta^0}_1-\norm{\widehat{\Psi}\beta^0+\widehat{\Psi}h_{T}}_1-\norm{\widehat{\Psi}h_{T^c}}_1\\
\label{boundgamma} &\le \norm{\widehat{\Psi}h_{T}}_1-\norm{\widehat{\Psi}h_{T^c}}_1.
\end{align}
We have an analogous bound for $\gamma^0$
\begin{equation} \norm{\gamma^0}_1- \norm{\widehat{\gamma}^0}_1 \le  \left|\left|f_{S}\right|\right|_1-\left|\left|f_{S^c}\right|\right|_1.
\label{boundalpha}
\end{equation}
The following holds 
\begin{align}\notag &(Q^0(\widehat{\beta}^0, \widehat{\gamma}^0))^{1/2} - (Q^0(\beta^0,\gamma^0))^{1/2}\\& 
\notag \le \frac{\lambda_{\beta}^0}{n}(\norm{\widehat{\Psi}\beta^0}_1-\norm{\widehat{\Psi}\widehat{\beta}^0}_1) + \frac{\lambda_{\gamma}^0}{n}(\norm{\gamma^0}_1-\norm{\widehat{\gamma}^0}_1)\\
\label{lbound} &\le \frac{\lambda_{\gamma}^0}{n}( \norm{\widehat{\Psi}h_{T}}_1-\norm{\widehat{\Psi}h_{T^c}}_1) +  \frac{\lambda_{\gamma}^0}{n}\left(\left|\left|f_{S}\right|\right|_1-\left|\left|f_{S^c}\right|\right|_1\right)
\end{align}
 By convexity, if $Q^0(\beta^0,\gamma^0)\ne 0$, it holds that
\begin{align}\notag(Q^0(\widehat{\beta}^0, \widehat{\gamma}^0))^{1/2} - (Q^0(\beta^0,\gamma^0))^{1/2}&\ge -\frac{1}{nQ^0(\beta^0,\gamma^0)^{1/2}}\left(\left<X^\top\xi^0,h \right>+\left<\xi^0,f \right>\right)\\
\label{rbound}&\ge  -\frac{1}{2n}\left(\lambda_{\beta}^0 \norm{\widehat{\Psi}h}_1+ \lambda_{\gamma}^0 \norm{f}_1\right).
\end{align}
This last inequality is also straightforwardly true when $Q^0(\beta^0,\gamma^0)=0$. Combining \eqref{lbound} and \eqref{rbound}, we get 
$$  -\frac{1}{2n}\left(\lambda_{\beta}^0 \norm{\widehat{\Psi}h}_1+ \lambda_{\gamma}^0 \left|\left|f\right|\right|_1\right)\le  \frac{\lambda_{\beta}^0}{n}( \norm{\widehat{\Psi}h_{T}}_1-\norm{\widehat{\Psi}h_{T^c}}_1) + \frac{\lambda_{\gamma}^0}{n}\left(\left|\left|f_{S}\right|\right|_1-\left|\left|f_{S^c}\right|\right|_1\right),$$
which implies that $(h,f) \in \mathbb{C}^0$.
Next, we have 
   \begin{align*}Q^0(\widehat{\beta}^0, \widehat{\gamma}^0)- Q^0(\beta^0,\gamma^0)&=\frac{1}{n}\norm{Xh + f}_2^2-\frac{2}{n}\left<\xi,Xh+f\right>\\
 &=\frac{1}{n}\norm{Xh + f}_2^2-\frac{2}{n}\left(\left<X^\top\xi^0,h \right>+\left<\xi^0,f \right>\right)\\
&\ge \frac{1}{n}\norm{Xh + f}_2^2-\frac{\norm{\xi^0}_2}{n^{3/2}}\left(\lambda_{\beta}^0 \norm{\widehat{\Psi}h}_1+ \lambda_{\gamma}^0 \left|\left|f\right|\right|_1\right)\\
&\ge\frac{1}{n}\norm{Xh + f}_2^2-\frac{4\norm{\xi^0}_2}{n^{3/2}}\left(\lambda_{\beta}^0 \norm{\widehat{\Psi}h_T}_1+ \lambda_{\gamma}^0 \left|\left|f_S\right|\right|_1\right).\end{align*}
Because $(h,f)\in\mathbb{C}^0$, \eqref{rbound} implies 
$(Q^0(\widehat{\beta}^0, \widehat{\gamma}^0))^{1/2} - (Q^0(\beta^0,\gamma^0))^{1/2}\ge -2n^{-1}\left(\lambda_{\beta}^0 \norm{\widehat{\Psi}h_T}_1+ \lambda_{\gamma}^0 \left|\left|f_S\right|\right|_1\right)$. Combined with \eqref{lbound}, it leads to 
$$\left|(Q^0(\widehat{\beta}^0, \widehat{\gamma}^0))^{1/2} - (Q^0(\beta^0,\gamma^0))^{1/2} \right|\le \frac2n \left(\lambda_{\beta}^0 \norm{\widehat{\Psi}h_T}_1+ \lambda_{\gamma}^0 \left|\left|f_S\right|\right|_1\right).$$
This yields
\begin{align*}
&Q^0(\widehat{\beta}^0, \widehat{\gamma}^0)- Q^0(\beta^0,\gamma^0)\\&\le \left((Q^0(\widehat{\beta}^0, \widehat{\gamma}^0))^{1/2} - (Q^0(\beta^0,\gamma^0))^{1/2}\right)  \left((Q^0(\widehat{\beta}^0, \widehat{\gamma}^0))^{1/2} + (Q^0(\beta^0,\gamma^0))^{1/2}\right)\\
&\le \frac2n\left(\lambda_{\beta^0} \norm{|\widehat{\Psi}h_T}_1+ \lambda_{\gamma}^0 \left|\left|f_S\right|\right|_1\right)  \left(2(Q^0(\beta^0,\alpha^0))^{1/2} +\frac2n\left(\lambda_{\beta}^0 \norm{\widehat{\Psi}h_T}_1+ \lambda_{\gamma}^0 \left|\left|f_S\right|\right|_1\right) \right)  
\end{align*}
Then, this implies 
   \begin{align*}&\frac{1}{n}\norm{Xh + f}_2^2\\
&\quad \le Q^0(\widehat{\beta}^0, \widehat{\gamma}^0)- Q^0(\beta^0,\gamma^0)+ \frac{4\norm{\xi^0}_2}{n^{3/2}}\left(\lambda_{\beta}^0 \norm{\widehat{\Psi}h_T}_1+ \lambda_{\gamma}^0 \left|\left|f_S\right|\right|_1\right)\\
 &\quad \le \frac4n\left(\lambda_{\beta}^0 \norm{\widehat{\Psi}h_T}_1+ \lambda_{\gamma}^0 \left|\left|f_S\right|\right|_1\right)^2 +\frac{8\norm{\xi}_2}{n^{3/2}}\left(\lambda_{\beta}^0 \norm{\widehat{\Psi}h_T}_1+ \lambda_{\gamma}^0 \left|\left|f_S\right|\right|_1\right).\\
&\quad \le \frac4n\left(\lambda_{\beta}^0 \sqrt{\norm{\beta^0}_0}\norm{\widehat{\Psi}h_T}_2+ \lambda_{\gamma}^0n\sqrt{\epsilon^k}\frac{1}{\sqrt{n}} \left|\left|f_S\right|\right|_2\right)^2 +\frac{8\norm{\xi}_2}{n^{3/2}}\left(\lambda_{\beta}^0 \sqrt{\norm{\beta^0}_0}\norm{\widehat{\Psi}h_T}_2+ \lambda_{\gamma}^0n\sqrt{\epsilon^k}\frac{1}{\sqrt{n}} \left|\left|f_S\right|\right|_2\right)\\
&\quad \le \left(\frac{2M^0}{n}\right)^2 \left(\norm{\widehat{\Psi}h}_2+\frac{1}{\sqrt{n}}\norm{f}_2\right)^2 + 8M^0\frac{\norm{\xi^0}_2}{n^{3/2}}\left(\norm{\widehat{\Psi}h}_2+\frac{1}{\sqrt{n}}\norm{f}_2\right)
\end{align*}
Because $(h,f)\in\mathbb{C}^0$, this implies 
$$(\kappa_*^0)^2\left(\norm{\widehat{\Psi}h}_2+\frac{1}{\sqrt{n}}\norm{f}_2\right)^2\le\left(\frac{2M^0}{n}\right)^2 \left(\norm{\widehat{\Psi}h}_2+\frac{1}{\sqrt{n}}\norm{f}_2\right)^2 + 8M^0\frac{\norm{\xi^0}_2}{n^{3/2}}\left(\norm{\widehat{\Psi}h}_2+\frac{1}{\sqrt{n}}\norm{f}_2\right).$$
If $\norm{\widehat{\Psi}h}_2+\norm{f}_2\ne 0$, we obtain
$$\left(\norm{\widehat{\Psi}h}_2+\frac{1}{\sqrt{n}}\norm{f}_2\right)\le\left((\kappa_*^0)^2-(2M^0/n)^2\right)^{-1}  8M^0\frac{\norm{\xi^0}_2}{n^{3/2}}.$$ 
Therefore, we have 
\begin{align*}
\norm{h}_1+\frac{\norm{f}_1}{\sqrt{n}}&\le \left( \sqrt{\norm{\beta^0}_0} + \sqrt{\epsilon^0} \right)(\lambda_{\max}(\widehat{\Psi})\vee 1)\left((\kappa_*^0)^2-(2M^0/n)^2\right)^{-1}  8M^0\frac{\norm{\xi^0}_2}{n^{3/2}}\\
\norm{h}_2+\frac{\norm{f}_2}{\sqrt{n}}&\le (\lambda_{\max}(\widehat{\Psi})\vee 1)\left((\kappa_*^0)^2-(2M^0/n)^2\right)^{-1}  8M^0\frac{\norm{\xi^0}_2}{n^{3/2}}\\
\norm{Xh+f}_2 &\le (\norm{X}_2\vee \sqrt{n})  (\lambda_{\max}(\widehat{\Psi})\vee 1)\left((\kappa_*^0)^2-(2M^0/n)^2\right)^{-1}  8M^0\frac{\norm{\xi^0}_2}{n^{3/2}}.
\end{align*}
To conclude, note that $\norm{\xi^0}_2= O_{\textrm{P}}(\sqrt{n})$ and $\Pr(\mathcal{E})\to 1$.

\subsection{Proof of Theorem \ref{than}}We work on the event 
$$\mathcal{E} =\left\{\lambda_{\beta}^0\ge 2\sqrt{n}\sup\limits_{j=1\dots,p}\frac{|(X^\top\xi^0)_j|}{\norm{\xi^0}_2\widehat{\Psi}_{jj}}\right\}\cap\left\{\lambda_{\gamma}^0\ge 2\sqrt{n}\frac{\norm{\xi^0}_\infty}{\norm{\xi^0}_2}\right\}$$
which has probability approaching 1 by \eqref{ascv} \label{cvi}. This implies that all the convergence in probability and distribution statements valid on this event will hold unconditional to this event as well. We use the notations $E=(\xi_1,\dots,\xi_n)^\top$ and $\widehat{E}=(\widehat{\xi}_1,\dots,\widehat{\xi}_n)^\top$.
\subsubsection{Proof of $\boldsymbol{\widehat{\Sigma}_\xi\xrightarrow{\Pr}\Sigma_\xi}$}

We have
\begin{align}
\notag \frac{1}{n}\widehat{E}^\top\widehat{E}&=\frac{1}{n}(\widehat{E}-E+E)^\top(\widehat{E}-E+E)\\
\label{Edecomposed}&= \frac{1}{n}\left[(\widehat{E}-E)^\top E +E^\top(\widehat{E}-E)+ (\widehat{E}-E)^\top(\widehat{E}-E)+E^\top E\right]
\end{align}
By Assumption \ref{asrate} and Theorem \ref{thcv}, it holds that 
\begin{equation}
 \norm{\widehat{E}-E}_2^2=\sum_{k=1}^K\norm{X(\widehat{\beta}^k-\beta^k)+\widehat{\gamma}^k-\gamma^k}_2^2
=O_\text{P}\left(\sum_{k=1}^K\frac{(M^k)^2}{n}\right)=o_\text{P}\left(\sqrt{n}\right)\label{cvE}
\end{equation}
Next, we have 
\begin{align*}
\left|( (\widehat{E}-E)^\top E)_{kk'}\right|&= \left|\left( X(\widehat{\beta}^k-\beta^k)+\widehat{\gamma}^k-\gamma^k\right)^\top\xi^{k'}\right|\\
&\le  \left|(\widehat{\beta}^k-\beta^k)^\top X^\top\xi^{k'}\right| +  \left|(\widehat{\gamma}^k-\gamma^k)^\top\xi^{k'}\right|\\
&\le \norm{X^\top\xi^{k'}}_\infty \norm{\widehat{\beta}^k-\beta^k}_1+ \norm{\xi^{k'}}_\infty \norm{\widehat{\gamma}^k-\gamma^k}_1\\
&\le \lambda^{k'}_\beta \frac{\norm{\xi^{k'}}_2}{2n} \norm{X}_{2,\infty} \norm{\widehat{\beta}^k-\beta^k}_1+ \lambda^{k'}_\gamma  \frac{\norm{\xi^{k'}}_2}{2\sqrt{n}} \norm{\widehat{\gamma}^k-\gamma^k}_1
\end{align*}
because we work on the event $\mathcal{E}$.
Therefore, by Theorem \ref{thcv}, it holds that
$$(\widehat{E}-E)^\top E=O_{\text{P}}\left(\sqrt{\bar \mu} \left(\frac{\norm{X}_{2,\infty}}{\sqrt{n}}\vee 1\right)\frac{\bar M}{n}\left(\bar\lambda_\beta\frac{\norm{X}_{2,\infty}}{\sqrt{n}} + \bar\lambda_\gamma \sqrt{n}\right)\right)=o_{\text{P}}\left(\sqrt{n}\right),$$
and similarly $E^\top(\widehat{E}-E)=o_{\textrm{P}}(\sqrt{n})$. Using \eqref{Edecomposed}, we get 
 \begin{equation}\label{Eestimate}
 \frac1n\widehat{E}^\top\widehat{E}- \frac1nE^\top E=o_{\textrm{P}}\left(\frac{1}{\sqrt{n}}\right)
\end{equation}
By the law of large numbers and Assumption \ref{asan}, we have $n^{-1}E^\top E \xrightarrow{\Pr}\Sigma_\xi$, which implies $\widehat{\Sigma}_\xi\xrightarrow{\Pr}\Sigma_\xi$ by \eqref{Eestimate}.
\subsubsection{Proof of asymptotic normality} 
 We have
\begin{align}
\notag \frac{1}{\sqrt{n}}\widehat{E}^\top\widehat{\xi}^0&=\frac{1}{\sqrt{n}}(\widehat{E}-E+E)^\top(\widehat{\xi}^0-\xi^0+\xi^0)\\
\label{xidecomposed}&= \frac{1}{\sqrt{n}}\left[(\widehat{E}-E)^\top \xi^0 +E^\top(\widehat{\xi}^0-\xi^0)+ (\widehat{E}-E)^\top(\widehat{\xi}^0-\xi^0)+E^\top \xi\right].
\end{align}
By Assumption \ref{asrate} and Theorem \ref{thcv}, it holds that \begin{equation}\norm{\widehat{\xi}^0-\xi^0}_2= \norm{X(\widehat{\beta}^0-\beta^0)+\widehat{\gamma}^0-\gamma^0}_2=O_\text{P}\left(\frac{M^0}{\sqrt{n}}\right)=o_{\text{P}}(n^{\frac14}),\label{cvxi}\end{equation}
which implies $|(\widehat{E}-E)^\top(\widehat{\xi}^0-\xi^0)|=o_{\textrm{P}}(\sqrt{n})$ by \eqref{cvE} and the inequality of Cauchy-Schwarz. 
Next, we have 
\begin{align*}
\left|( (\widehat{\xi}^0-\xi^0)^\top E)_{k}\right|&= \left|\left( X(\widehat{\beta}^0-\beta^0)+\widehat{\gamma}^0-\gamma^0\right)^\top\xi^{k}\right|\\
&\le  \left|(\widehat{\beta}^0-\beta^0)^\top X^\top\xi^{k}\right| +  \left|(\widehat{\gamma}^0-\gamma^0)^\top\xi^{k}\right|\\
&\le \norm{X^\top\xi^{k}}_\infty \norm{\widehat{\beta}^0-\beta^0}_1+ \norm{\xi^{k}}_\infty \norm{\widehat{\gamma}^0-\gamma^0}_1\\
&\le \lambda^{k}_\beta \frac{\norm{\xi^{k}}_2}{2n}\norm{X}_{2,\infty} \norm{\widehat{\beta}^0-\beta^0}_1+ \lambda^{k}_\gamma  \frac{\norm{\xi^{0}}_2}{2\sqrt{n}} \norm{\widehat{\gamma}^k-\gamma^k}_1.
\end{align*}
Therefore, by Theorem \ref{thcv}, it holds that
$$(\widehat{\xi}^0-\xi^0)^\top E=O_{\text{P}}\left(\sqrt{\bar \mu}\left(\frac{\norm{X}_{2,\infty}}{\sqrt{n}}\vee 1\right) \frac{\bar M}{n}\left(\bar\lambda_\beta\frac{\norm{X}_{2,\infty}}{\sqrt{n}} + \bar\lambda_\gamma   \sqrt{n}\right)\right)=o_{\text{P}}\left(\sqrt{n}\right).$$
Similarly, it holds that 
$$(\widehat{E}-E)^\top \xi^0=O_{\text{P}}\left(\sqrt{\bar \mu} \left(\frac{\norm{X}_{2,\infty}}{\sqrt{n}}\vee 1\right)\frac{\bar M}{n}\left(\bar\lambda_\beta\frac{\norm{X}_{2,\infty}}{\sqrt{n}} + \bar\lambda_\gamma  \sqrt{n}\right)\right)=o_{\text{P}}\left(\sqrt{n}\right).$$
Using \eqref{xidecomposed}, we get 
 \begin{equation}\label{xiestimate}
 \frac{1}{\sqrt{n}}\widehat{E}^\top\widehat{\xi}^0-  \frac{1}{\sqrt{n}}E^\top \xi^0=o_{\textrm{P}}\left(1\right)
\end{equation}
By \eqref{Eestimate}, the law of large numbers and the continuous mapping theorem, it holds that $n(\widehat{E}^\top\widehat{E})^{-1} \xrightarrow{\Pr} \Sigma_\xi^{-1}$ and $n(E^\top E )^{-1} \xrightarrow{\Pr} \Sigma_\xi^{-1}$, which implies $n(\widehat{E}^\top\widehat{E})^{-1} -n(E^\top E)^{-1}= o_{\textrm{P}}(1)$. This and \eqref{xiestimate} yield
\begin{align*}
\sqrt{n}\widehat{\alpha}&= \left(\frac1n\widehat{E}^\top\widehat{E}\right)^{-1} \frac{1}{\sqrt{n}}\widehat{E}^\top\widehat{\xi}^0 \\
&=  \left(\frac1n\widehat{E}^\top\widehat{E}\right)^{-1}\left( \frac{1}{\sqrt{n}}\widehat{E}^\top\widehat{\xi}^0-\frac{1}{\sqrt{n}}E^\top\xi^0\right)+\left( \left(\frac1n\widehat{E}^\top\widehat{E}\right)^{-1}-\left(\frac1n E^\top  E\right)^{-1} \right)\frac{1}{\sqrt{n}}E^\top \xi^0\\
&\quad+\left(\frac1n E^\top E\right)^{-1} \frac{1}{\sqrt{n}}E^\top \xi^0\\
&=o_{\textrm{P}}(1)  + \sqrt{n} \alpha+\left(\frac1n E^\top E\right)^{-1} \frac{1}{\sqrt{n}}E^\top u.
\end{align*}
We conclude using the central limit theorem and Slutsky's theorem.

\subsubsection{Proof of $\boldsymbol{\widehat{\sigma}\xrightarrow{\Pr}\sigma}$}
Let $\widehat{u}_i= \widehat{\xi}^0_i- \sum_{k=1}^K\widehat{\alpha}_k\widehat{\xi}_i^k$ and $\widehat{u}=(\widehat{u}_1,\dots, \widehat{u}_n)^\top$. We have 
\begin{align*}
\widehat{u}-u&= \widehat{\xi}^0-\xi^0- \sum_{k=1}^K(\widehat{\alpha}_k\widehat{\xi}^k-\alpha_k\xi^k)\\
&=  \widehat{\xi}^0-\xi^0- \sum_{k=1}^K(\widehat{\alpha}_k-\alpha _k)(\widehat{\xi}^k-\xi^k) + (\widehat{\alpha}_k-\alpha _k)\xi^k + \alpha_k (\widehat{\xi}^k-\xi^k)= o_{\textrm{P}}(n^{\frac14}),
\end{align*}
by \eqref{cvE}, \eqref{cvxi} and the fact that $\widehat{\alpha}-\alpha=o_{\textrm{P}}(1)$. This implies that \begin{align*}\widehat{\sigma}^2&=\frac1n\widehat{u}^\top \widehat{u}\\
&=\frac1n(\widehat{u}-u)^\top u+\frac1n u^\top(\widehat{u}-u)^\top+\frac1n(\widehat{u}-u)^\top(\widehat{u}-u)+ \frac1n u^\top u+\\
&= \frac1n u^\top u +o_{\textrm{P}}(1)\xrightarrow{\Pr} \sigma^2\end{align*}
by the law of large numbers.
\bibliographystyle{plainnat}
\bibliography{outliers_high_dim}

 \end{document}